\documentclass{article}
\usepackage{amscd}
\usepackage{amsmath}
\usepackage{amssymb}
\usepackage{amsthm}
\usepackage{ascmac}
\usepackage{bm}
\usepackage{graphicx}
\usepackage{here}
\usepackage[pdftex]{hyperref}%
\usepackage{mathtools}\numberwithin{equation}{section}
\usepackage{minitoc}
\usepackage{multicol}
\usepackage{multirow}
\usepackage{setspace}
\usepackage[dvipdfmx]{color}
\theoremstyle{definition}
\newtheorem{df}{Definition}
\newtheorem{note}{Notation}
\newtheorem{q}{Problem}
\newtheorem{eg}{Example}
\newtheorem{rmk}{Remark}
\theoremstyle{plain}
\newtheorem{thm}{Theorem}
\newtheorem{prop}{Proposition}

\newtheorem{cor}{Corollary}
\newenvironment{pf}{\begin{proof}}{\end{proof}}
\usepackage{tikz}
\usetikzlibrary{calc,intersections,patterns,through}
\newcommand{\ang}{\tikz@ang}
\def\tikz@ang(#1)(#2)#3{%
\pgfmathanglebetweenpoints{%
\pgfpointanchor{#1}{center}}{%
\pgfpointanchor{#2}{center}}
\pgfmathsetmacro{#3}{\pgfmathresult}%
}

\title{Rational Angle Bisection Problem in Higher Dimensional Spaces and Incenters of Simplices over Fields}
\author{Takashi HIROTSU}
\date{\today}
\begin{document}
\maketitle
\begin{abstract}
In this article, we generalize the following problem, which is called the rational angle bisection problem, to the $n$-dimensional space $k^n$ over a subfield $k$ of $\mathbb R$: in the coordinate plane, for which rational numbers $a$ and $b$ are the slopes of the angle bisectors between the two lines with slopes $a$ and $b$ rational? 
First, we provide several characterizations of when the angle bisectors between two lines with direction vectors in $k^n$ have direction vectors in $k^n.$ 
To find solutions to the problem in the case when $k = \mathbb Q,$ we derive a formula for the integral solutions of $x_1{}^2+\dots +x_n{}^2 = dx_{n+1}{}^2,$ which is a generalization of negative Pell's equation $x^2-dy^2 = -1,$ where $d$ is a square-free positive integer. 
Second, by applying the above characterizations, we establish a necessary and sufficient condition for the incenter of a given $n$-simplex with $k$-rational vertices to be $k$-rational. 
In the coordinate plane, we prove that every triangle with $k$-rational vertices and incenter can be obtained by scaling a triangle with $k$-rational side lengths and area, which is a generalization of a Heronian triangle. 
We also discuss certain fundamental properties of a few centers of a given triangle with $k$-rational vertices.
\end{abstract}
\section{Introduction}\label{sec-intro}
Throughout this article, let $n \geq 2$ be an integer, and let $k$ be a subfield of $\mathbb R.$ 
We consider the following problem in $\mathbb R^n,$ which is a generalization of the rational angle bisection problem in $\mathbb R^2$ (see \cite{Hir24}).
\begin{q}\label{q-rat-bisec}
For which $\bm{a},$ $\bm{b} \in k^n$ do the angle bisectors between the two lines with direction vectors $\bm{a}$ and $\bm{b}$ have direction vectors in $k^n$?
\end{q}
We denote the inner product of $\bm{a},$ $\bm{b} \in k^n$ by $\langle\bm{a},\bm{b}\rangle.$
\begin{rmk}
Given two lines, we consider the two angles formed by them, regardless of whether they are acute or not. 
If the bisector of one of the angles between $\bm{a},$ $\bm{b} \in k^n$ has a direction vector $\bm{c} \in k^n,$ then that of the other is perpendicular to $\bm{c}$ and has direction vector $\langle\bm{b},\bm{c}\rangle\bm{a}-\langle\bm{a},\bm{c}\rangle\bm{b} \in k^n.$
\end{rmk}
In the case when $k/\mathbb Q$ is algebraic, let $O_k$ be the integer ring of $k.$ 
In this case, Problem \ref{q-rat-bisec} is reduced to the following problem, which becomes an important problem involving solutions of equations that generalize Pythagorean and negative Pell's equations as discussed later.
\begin{q}\label{q-int-bisec}
For which $\bm{a},$ $\bm{b} \in O_k{}^n$ do the angle bisectors between the two lines with direction vectors $\bm{a}$ and $\bm{b}$ have direction vectors in $O_k{}^n$?
\end{q}
Problem \ref{q-int-bisec} is important in drawing techniques when $n = 2$ and $k = \mathbb Q,$ and has an application to optical engineering when $n = 3$ and $k = \mathbb Q$; that is, we can specify the radiation range and the axis of light with a ratio of integers without the errors arising from approximations of irrational numbers by using a solution to Problem \ref{q-int-bisec}.\par
Let $k^{\times 2} = \{ q^2 \mid q \in k^\times\}.$ 
The following theorem plays an important role in solving Problems \ref{q-rat-bisec} and \ref{q-int-bisec}.
\begin{thm}\label{thm-ang}
Let $\bm{a}$ and $\bm{b}$ be linearly independent vectors of $k^n.$
\begin{enumerate}
\item[{\rm (1)}]
The following conditions are equivalent. 
\begin{enumerate}
\item[{\rm (A1)}]
One of the angle bisectors between the two lines with direction vectors $\bm{a}$ and $\bm{b}$ has direction vector $\bm{c} \in k^n.$
\item[{\rm (A2)}]
One of the angle bisectors between the two hyperplanes with normal vectors $\bm{a}$ and $\bm{b}$ has normal vector $\bm{c} \in k^n.$
\item[{\rm (A3)}]
There exists a vector $\bm{c} \in k^n$ such that 
\begin{equation} 
\langle\bm{a},\bm{c}\rangle ^2|\bm{b}|^2 = \langle\bm{b},\bm{c}\rangle ^2|\bm{a}|^2. \label{eq-star} 
\end{equation}
\item[{\rm (A4)}]
We have 
\[ |\bm{a}|^2 \equiv |\bm{b}|^2 \pmod{k^{\times 2}}.\]
\end{enumerate}
If {\rm (A4)} holds, then $\bm{c}$ in {\rm (A1)}--{\rm (A3)} can be taken as 
\begin{equation} 
\bm{c} = \bm{a}\pm\frac{|\bm{a}|}{|\bm{b}|}\bm{b}. \label{eq-c} 
\end{equation}
\item[{\rm (2)}]
Suppose that $k/\mathbb Q$ is algebraic and $\bm{a} = (a_1,\dots,a_n),$ $\bm{b} = (b_1,\dots,b_n) \in O_k{}^n.$ 
Then the following conditions are equivalent. 
\begin{enumerate}
\item[{\rm (B1)}]
One of the angle bisectors between the two lines with direction vectors $\bm{a}$ and $\bm{b}$ has direction vector $\bm{c} \in O_k{}^n.$
\item[{\rm (B2)}]
One of the angle bisectors between the two hyperplanes with normal vectors $\bm{a}$ and $\bm{b}$ has normal vector $\bm{c} \in O_k{}^n.$
\item[{\rm (B3)}]
There exists a vector $\bm{c} \in O_k{}^n$ such that \eqref{eq-star}. 
\item[{\rm (B4)}]
There exist algebraic integers $a_{n+1},$ $b_{n+1},$ $d \in O_k$ such that $d > 0$ and $(x_1,\dots,x_n,x_{n+1}) =$\\
$(a_1,\dots,a_n,a_{n+1}),$ $(b_1,\dots,b_n,b_{n+1})$ are solutions of 
\begin{equation} 
\sum_{i = 1}^nx_i{}^2 = dx_{n+1}{}^2. \label{eq-pell} 
\end{equation} 
\end{enumerate}
If {\rm (B4)} holds, then $\bm{c}$ in {\rm (B1)}--{\rm (B3)} can be taken as 
\[\bm{c} = b_{n+1}\bm{a}\pm a_{n+1}\bm{b}.\] 
\end{enumerate}
\end{thm}
We prove Theorem \ref{thm-ang} in Section \ref{sec-ang}.
\begin{eg}
Let $\bm{a} = (a_1,\dots,a_n),$ $\bm{b} = (b_1,\dots,b_n) \in k^n.$ 
If there exists a nonnegative number $q \in k$ such that 
\[\{ |b_1|,\dots,|b_n|\} = \{ q|a_1|,\dots,q|a_n|\} \quad\text{or}\quad \{ |a_1|,\dots,|a_n|\} = \{ q|b_1|,\dots,q|b_n|\},\] 
then $\bm{a},$ $\bm{b},$ and $\bm{c} = \bm{a}$ satisfy \eqref{eq-star}. 
We call such a solution of \eqref{eq-star} {\it trivial}.
\end{eg}
Several examples of solutions to Problem \ref{q-int-bisec} in the case when $n = 2$ and $k = \mathbb Q$ are given in \cite{Hir24}. 
Here, we provide a few solutions to Problem \ref{q-int-bisec} derived from nontrivial solutions of \eqref{eq-star} in the case when $n = 3$ and $k = \mathbb Q.$
\begin{eg}
\begin{enumerate}
\item[(1)]
If $\bm{a} = (2,2,1)$ and $\bm{b} = (3,6,2),$ then the angle bisectors between the two lines with direction vectors $\bm{a}$ and $\bm{b}$ have direction vectors $7\bm{a}\pm 3\bm{b},$ since $|\bm{a}| = 3$ and $|\bm{b}| = 7$ imply that $7\bm{a}$ and $3\bm{b}$ span a rhomb (see Figure 1).
\begin{figure}[h]\label{fig-rhomb}
\centering
\includegraphics{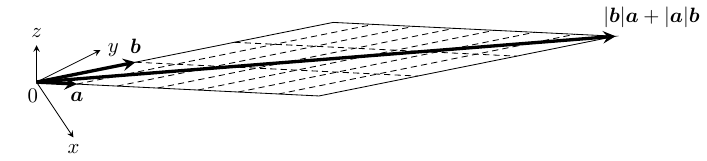}
\caption{The bisector of the acute angle between $\bm{a} = (2,2,1)$ and $\bm{b} = (3,6,2).$}
\end{figure}
\item[(2)]
If $\bm{a} = (1,1,4)$ and $\bm{b} = (3,4,5),$ then the angle bisectors between the two lines with direction vectors $\bm{a}$ and $\bm{b}$ have direction vectors $5\bm{a}\pm 3\bm{b},$ since $|\bm{a}| = 3\sqrt 2$ and $|\bm{b}| = 5\sqrt 2$ imply that $5\bm{a}$ and $3\bm{b}$ span a rhomb.
\end{enumerate}
\end{eg}\par
To provide more solutions to Problem \ref{q-int-bisec} in the case when $k = \mathbb Q,$ we derive a formula for \eqref{eq-pell} in Section \ref{sec-plike}.\par
Next, we apply Theorem \ref{thm-ang} to the geometry of simplices. 
We use the following definition.
\begin{df}
\begin{enumerate}
\item[(1)]
A point $P \in \mathbb R^n$ is said to be {\it $k$-rational} if $P \in k^n.$
\item[(2)]
A polytope is said to be {\it over $k$} if all of its vertices are $k$-rational. 
\item[(3)]
The volume $V$ of a polytope is said to be {\it $k$-rational} if $V \in k.$
\item[(4)]
An $n$-simplex is said to be {\it Heronian over $k$} if the volumes of its $i$-faces are all $k$-rational for each $i \in \{ 1,\dots,n\}.$
\end{enumerate}
\end{df}
\begin{thm}\label{thm-sim}
In $\mathbb R^n,$ given $n$-simplex $A_0A_1\cdots A_n$ over $k,$ let $a_i$ be the volume of the facet opposite to $A_i$ for each $i \in \{ 0,1,\dots,n\}.$ 
The following conditions are equivalent.
\begin{enumerate}
\item[{\rm (S1)}]
The incenter $I$ of $n$-simplex $A_0A_1\cdots A_n$ is $k$-rational.
\item[{\rm (S2)}]
We have 
\begin{equation} 
a_0{}^2 \equiv a_1{}^2 \equiv \cdots \equiv a_n{}^2 \pmod{k^{\times 2}}. \label{eq-rat-sim}
\end{equation} 
\end{enumerate}
If $k/\mathbb Q$ is algebraic, then these conditions are equivalent to the following condition.
\begin{enumerate}
\item[{\rm (S2)${}'$}]
There exists a positive algebraic integer $d \in O_k$ such that 
\begin{equation} 
a_0{}^2 \equiv a_1{}^2 \equiv \cdots \equiv a_n{}^2 \equiv d\pmod{k^{\times 2}}.\label{eq-int-sim}
\end{equation} 
\end{enumerate}
\end{thm}
\begin{rmk}
In \eqref{eq-rat-sim} or \eqref{eq-int-sim}, we have $a_0{}^2,$ $a_1{}^2,$ $\dots,$ $a_n{}^2 \in k^\times$ by the Cayley--Menger formula (see \cite[Section 5, Chapter VIII]{Som58}).
\end{rmk}
\begin{thm}\label{thm-trig}
In $\mathbb R^2,$ given triangle $ABC$ over $k$ with side lengths $BC = a,$ $CA = b,$ and $AB = c,$ inradius $r,$ and exradii $r_A,$ $r_B,$ and $r_C,$ the following conditions are equivalent.
\begin{enumerate}
\item[{\rm (T1)}]
At least one of the incenter $I,$ the excenters $I_A,$ $I_B,$ and $I_C$ relative to $A,$ $B,$ and $C,$ the Gergonne point $Ge,$ and the Nagel point $Na$ of triangle $ABC$ is $k$-rational.
\item[{\rm (T2)}]
We have 
\begin{equation} 
a^2 \equiv b^2 \equiv c^2 \pmod{k^{\times 2}}. \label{eq-trig} 
\end{equation} 
\end{enumerate}
If these equivalent conditions hold, then all of $I,$ $I_A,$ $I_B,$ $I_C,$ $Ge,$ and $Na$ are $k$-rational, and triangle $ABC$ can be obtained by scaling a Heronian triangle over $k$ by a factor $q\sqrt d$ for some positive numbers $d,$ $q \in k,$ and $a^2,$ $b^2,$ $c^2,$ $r^2,$ $r_A{}^2,$ $r_B{}^2,$ and $r_C{}^2$ are numbers in $k^\times$ equivalent to $d$ modulo $k^{\times 2}.$
\end{thm}
We prove Theorem \ref{thm-sim} in Section \ref{sec-sim}.  
We also prove Theorem \ref{thm-trig} and provide a few examples of triangles with $k$-rational vertices and incenter in Section \ref{sec-trig}.
\section{Angle Bisectors}\label{sec-ang}
In this section, we prove Theorem \ref{thm-ang}.
\begin{pf}[Proof of Theorem \ref{thm-ang}]
We prove (1).\par
First, (A1) and (A2) are equivalent, since the angle between two hyperplanes is defined as the smaller of the angle between their normal vectors and its supplementary angle.\par
Next, we prove that (A1) and (A3) are equivalent. 
Let $\bm{a} = \overrightarrow{OA_+},$ $\bm{b} = \overrightarrow{OB},$ and $\bm{c} = \overrightarrow{OC}$ with $O \in \mathbb R^n$ and $A_+,$ $B,$ $C \in \mathbb R^n\setminus\{ O\}.$ 
Let $A_-$ be the point symmetric to $A_+$ with respect to $O.$ 
The condition that $\angle A_+OC = \angle BOC$ or $\angle A_-OC = \angle BOC$ is equivalent to 
\begin{align*} 
\frac{\langle\pm\bm{a},\bm{c}\rangle}{|\pm\bm{a}||\bm{c}|} &= \frac{\langle\bm{b},\bm{c}\rangle}{|\bm{b}||\bm{c}|}, \\ 
\pm\frac{\langle\bm{a},\bm{c}\rangle}{|\bm{a}||\bm{c}|} &= \frac{\langle\bm{b},\bm{c}\rangle}{|\bm{b}||\bm{c}|}, 
\end{align*} 
and \eqref{eq-star} by multiplying both sides by $|\bm{a}||\bm{b}||\bm{c}|$ and squaring both sides.\par
Furthermore, it is obvious that (A3) implies (A4). 
Conversely, if (A4) holds, then (A3) holds since the vector $\bm{c} \in k^n$ defined by \eqref{eq-c} satisfies \eqref{eq-star}.\par 
We prove (2). 
By a similar argument as above, we can clarify that (B1)--(B3) are equivalent. 
Suppose that there exists a vector $\bm{c} \in O_k{}^n$ such that \eqref{eq-star}. 
Then we have $|\bm{a}|^2 \equiv |\bm{b}|^2 \pmod{k^{\times 2}},$ which implies that there exist algebraic integers $a_{n+1},$ $b_{n+1},$ $d \in O_k$ such that $d > 0$ and 
\begin{equation} 
|\bm{a}|^2 = da_{n+1}{}^2 \quad\text{and}\quad |\bm{b}|^2 = db_{n+1}{}^2. \label{eq-pels} 
\end{equation} 
This implies that $(x_1,\dots,x_n,x_{n+1}) = (a_1,\dots,a_n,a_{n+1}),$ $(b_1,\dots,b_n,b_{n+1})$ are solutions of \eqref{eq-pell}. 
The converse is also true.\qedhere
\end{pf}
\section{Pythagorean and Pell-like Equations}\label{sec-plike}
In this section, we consider Equation \eqref{eq-pell} in the case when $k = \mathbb Q.$ 
We use the following symbols.
\begin{note}
Let $\mathbb N$ denote the additive monoid of all nonnegative integers. 
For each prime number $p,$ let $\mathrm{ord}_p:\mathbb Q^\times\to\mathbb Z$ be the normalized $p$-adic additive valuation. 
We denote the greatest common divisor of $a,$ $b \in \mathbb Z\setminus\{ 0\}$ by $\mathrm{gcd}\,(a,b).$ 
Let $\mathfrak S_n$ denote the symmetric group of degree $n.$
\end{note}
If $d = 1,$ then \eqref{eq-pell} is a generalization of the Pythagorean equation $x^2+y^2 = z^2,$ and a positive integral solution of \eqref{eq-pell} is called a {\it Pythagorean $(n+1)$-tuple}. 
The following formula is well-known. 
\begin{thm}[{\cite[Section 10]{Car15}}]\label{thm-form-py}
Every integral solution $(x_1,x_2,x_3,x_4)$ of $x_1{}^2+x_2{}^2+x_3{}^2 = x_4{}^2$ is given by 
\begin{equation} 
(x_{\sigma (1)},x_{\sigma (2)},x_{\sigma (3)},x_4) = (s(t^2+u^2-v^2-w^2),s\cdot 2(tv+uw),s\cdot 2(tw-uv),\pm s(t^2+u^2+v^2+w^2)) \label{eq-form-py} 
\end{equation} 
for some $\sigma \in \mathfrak S_3$ and $(s,t,u,v,w) \in \mathbb Z^5.$
\end{thm}
\begin{rmk}
Euclid's well-known formula for Pythagorean triples is obtained by letting $s > 0,$ $t > v > u = w = 0$ and $x_4 > 0$ in \eqref{eq-form-py}.
\end{rmk}
Combining Theorem \ref{thm-ang} with Theorem \ref{thm-form-py}, we obtain the following solutions to Problem \ref{q-int-bisec}.
\begin{eg}
If $(a_{\sigma (1)},a_{\sigma (2)},a_{\sigma (3)}) = (1,0,0)$ and $(b_{\tau (1)},b_{\tau (2)},b_{\tau (3)}) = (t^2+u^2-v^2-w^2,2(tv+uw),2(tw-uv)),$ where $\sigma,$ $\tau \in \mathfrak S_3$ and $(t,u,v,w) \in \mathbb Z^4\setminus \{ (0,0,0,0)\},$ then $\bm{a} = (a_1,a_2,a_3),$ $\bm{b} = (b_1,b_2,b_3),$ and $\bm{c} = (t^2+u^2+v^2+w^2)\bm{a}\pm\bm{b}$ satisfy \eqref{eq-star}.
\end{eg}
\begin{eg}
If $(x_1,x_2,x_3,x_4) = (a_1,a_2,a_3,a_4),$ $(b_1,b_2,b_3,b_4)$ are integral solutions of $x_1{}^2+x_2{}^2+x_3{}^2 = x_4{}^2$ for some $a_4,$ $b_4 \in \mathbb Z\setminus\{ 0\},$ then $\bm{a} = (a_1,a_2,a_3),$ $\bm{b} = (b_1,b_2,b_3),$ and $\bm{c} = (b_4\bm{a}\pm a_4\bm{b})/g$ satisfy \eqref{eq-star}, where $g = \mathrm{gcd}(a_4,b_4).$
We summarize a few examples in the following table. 
\begin{center}\begin{tabular}{|c||c|c|c|c|c|c|c|c|c|} \hline 
$\bm{a}$ & $(3,4,0)$ & $(3,4,0)$ & $(63,16,0)$ & $(3,4,0)$ & $(3,4,0)$ & $(5,12,0)$ & $(1,2,2)$ & $(1,2,2)$ & $(1,8,4)$ \\ 
$a_4$ & $5$ & $5$ & $65$ & $5$ & $5$ & $13$ & $3$ & $3$ & $9$ \\ 
$\bm{b}$ & $(5,0,12)$ & $(7,0,24)$ & $(33,0,56)$ & $(1,2,2)$ & $(5,2,14)$ & $(3,4,12)$ & $(3,2,6)$ & $(1,8,4)$ & $(7,4,4)$ \\ 
$b_4$ & $13$ & $25$ & $65$ & $3$ & $15$ & $13$ & $7$ & $9$ & $9$ \\ 
$g$ & $1$ & $5$ & $65$ & $1$ & $5$ & $13$ & $1$ & $3$ & $9$ \\ \hline 
\end{tabular}\end{center}
\end{eg}
If $d > 1,$ then \eqref{eq-pell} is equivalent to $x_n{}^2-dx_{n+1}{}^2 = -\sum_{i = 1}^{n-1}x_i{}^2,$ which is a generalization of negative Pell's equation $x^2-dy^2 = -1$ and is called a {\it Pell-like equation} in this article. 
The following theorem follows from Fermat's two-square theorem, Legendre's three-square theorem, and Lagrange's four-square theorem. 
\begin{thm}\label{thm-ex-pell}
Let $d > 1$ be a square-free integer.
\begin{enumerate}
\item[{\rm (1)}]
Suppose that $n = 2.$ 
Then \eqref{eq-pell} has an integral solution if and only if $d$ has no prime divisors congruent to $3$ modulo $4.$
\item[{\rm (2)}]
Suppose that $n = 3.$ 
Then \eqref{eq-pell} has an integral solution if and only if $d \not\equiv 7 \pmod 8.$
\item[{\rm (3)}]
Suppose that $n \geq 4.$ 
Then \eqref{eq-pell} has an integral solution independently of the value of $d.$ 
In particular, if $n \geq 5,$ then \eqref{eq-pell} has an integral solution $(x_1,\dots,x_n,x_{n+1})$ such that $(x_{\sigma (1)},\dots,x_{\sigma (n-4)}) = (a_1,\dots,a_{n-4})$ for any $\sigma \in \mathfrak S_n$ and $(a_1,\dots,a_{n-4}) \in \mathbb Z^{n-4}.$
\end{enumerate}
\end{thm}
\begin{pf}
To prove (1), suppose that $n = 2.$ 
If $d$ has no prime divisors congruent to $3$ modulo $4,$ then \eqref{eq-pell} has an integral solution with $x_3 = 1$ by Fermat's two-square theorem. 
Conversely, if \eqref{eq-pell} has an integral solution $(x_1,x_2,x_3),$ then every prime divisor $p$ of $dx_3{}^2$ such that $p \equiv 3 \pmod 4$ satisfies $\mathrm{ord}_p(dx_3{}^2) \equiv 0 \pmod 2$ by Fermat's two-square theorem, which implies $\mathrm{ord}_p(d) = 0.$\par
To prove (2), suppose that $n = 3.$ 
If $d \not\equiv 7 \pmod 8,$ then \eqref{eq-pell} has an integral solution with $x_4 = 1$ by Legendre's three-square theorem. 
Conversely, if \eqref{eq-pell} has an integral solution $(x_1,x_2,x_3,x_4),$ then $dx_4{}^2$ can not be written in the form $4^a(8b+7)$ with $a,$ $b \in \mathbb N$ by Legendre's three-square theorem, which implies $d \not\equiv 7 \pmod 8$ since $x_4{}^2/4^i = 4j(j+1)+1 \equiv 1 \pmod 8$ if $x_4 = 2^i(2j+1),$ $i \in \mathbb N,$ and $j \in \mathbb Z.$\par
Finally, (3) follows from Lagrange's four-square theorem. 
In particular, the latter half is true, since $x_{\sigma (n-3)}{}^2+x_{\sigma (n-2)}{}^2+x_{\sigma (n-1)}{}^2+x_{\sigma (n)}{}^2 = da_{n+1}{}^2-\sum_{i = 1}^{n-4}a_i{}^2$ has an integral solution for some $a_{n+1} \in \mathbb N$ such that $da_{n+1}{}^2 \geq \sum_{i = 1}^{n-4}a_i{}^2.$
\end{pf}
To provide solutions to Problem \ref{q-int-bisec}, we prove the theorem below, which is a generalization of the known result for the case when $n = 2$ (see \cite[Theorem 2]{Hir24}). 
To state the theorem, we use the following definition and symbols.
\begin{df}[{\cite[Definition 1]{Hir24}}]
Let $d > 1$ be a square-free integer, and let $z > 0$ be an integer.
\begin{enumerate}
\item[(1)]
Let $(x,y)$ be an integral solution of $|x^2-dy^2| = z.$ 
We say that $(x,y)$ is {\itshape strictly primitive}, if $x$ and $dy$ are coprime.
\item[(2)]
Let $(x,y) = (a_1,a_2),$ $(b_1,b_2)$ be positive integral solutions of $|x^2-dy^2| = z.$ 
We say that $(a_1,a_2)$ is {\it smaller} than $(b_1,b_2),$ if $a_2 < b_2,$ or if $a_1 < b_1$ and $a_2 = b_2.$
\item[(3)]
For each equation of $x^2-dy^2 = z,$ $x^2-dy^2 = -z,$ and $|x^2-dy^2| = z,$ we call its minimum positive integral solution its {\itshape fundamental solution}.
\end{enumerate}
\end{df}
\begin{note}[{\cite[Notation 2]{Hir24}}]
Let $d > 1$ be a square-free integer. 
Let $\eta$ be the fundamental unit of the real quadratic field $\mathbb Q(\sqrt d).$ 
We denote the $m$-th smallest positive integral solution of $|x^2-dy^2| = 1$ by $(x,y) = (f_m^{(d)},g_m^{(d)}).$ 
Let $S(d)$ denote the set of every prime number $p$ such that $|x^2-dy^2| = p^l$ has a strictly primitive integral solution for some integer $l > 0.$ 
For each $p \in S(d),$ let $l_p$ be the minimum integer $l > 0$ such that $|x^2-dy^2| = p^l$ has a strictly primitive integral solution, and let 
\[\xi _p = x_p+y_p\sqrt d\] 
with the fundamental solution $(x,y) = (x_p,y_p)$ of 
\[\begin{cases} 
x^2-dy^2 = p^{l_p} & \text{if }x^2-dy^2 = -1\text{ has an integral solution,} \\ 
|x^2-dy^2| = p^{l_p} & \text{otherwise}. 
\end{cases}\] 
Let $S(d)_-$ denote the set of every prime number $p$ such that $p \in S(d)$ and $x_p{}^2-dy_p{}^2 = -p^{l_p}.$ 
Furthermore, for each $\alpha = a_1+a_2\sqrt d \in \mathbb Q(\sqrt d),$ where $a_1,$ $a_2 \in \mathbb Q,$ we denote the conjugate and norm of $\alpha$ by 
\[\alpha ' = a_1-a_2\sqrt d \quad\text{and}\quad  N(\alpha ) = \alpha\alpha ' = a_1{}^2-da_2{}^2,\] 
respectively.
\end{note}
\begin{rmk}
The following fact is known (see \cite[Theorem 8]{Hir24}).
\begin{enumerate}
\item[\textup{(1)}]
If $d \equiv 1\ (\mathrm{mod}\ 8),$ or $d \equiv 5\ (\mathrm{mod}\ 8)$ and $\eta \in \mathbb Z[\sqrt d],$ or $d \equiv 2,$ $3\ (\mathrm{mod}\ 4),$ then $S(d)$ consists of all prime numbers that split in $\mathbb Q(\sqrt d).$
\item[\textup{(2)}]
If $d \equiv 5\ (\mathrm{mod}\ 8)$ and $\eta \notin \mathbb Z[\sqrt d],$ then $S(d)$ consists of all prime numbers that split in $\mathbb Q(\sqrt d)$ and $2.$
\end{enumerate}
\end{rmk}
\begin{thm}\label{thm-form-plike}
Let $d > 1$ be a square-free integer. 
Let $\sigma \in \mathfrak S_n,$ $(a_1,\dots,a_{n-1}) \in \mathbb Z^{n-1},$ and $z = \sum_{i = 1}^{n-1}a_i{}^2.$ 
If we have 
\[\mathrm{ord}_p(z) \equiv 0 \pmod 2\] 
for each prime number $p \notin S(d)$ and there exists an integer $m_p \geq 0$ such that 
\begin{equation} 
l_pm_p \equiv \mathrm{ord}_p(z) \pmod 2 \quad\text{and}\quad m_p \leq \mathrm{ord}_p(z)/l_p \label{eq-exp-cond} 
\end{equation} 
for each $p \in S(d),$ then \eqref{eq-pell} has an integral solution $(x_1,\dots,x_n,x_{n+1})$ such that $(x_{\sigma (1)},\dots,x_{\sigma (n-1)}) = (a_1,\dots,a_{n-1})$ and $\mathrm{gcd}\,(x_{\sigma (n)},d) = 1,$ and any such solution satisfies 
\begin{equation} 
x_{\sigma (n)}+x_{n+1}\sqrt d = \pm\eta ^m\prod_{p \in S(d)}\xi _p^*{}^{m_p}p^{(\mathrm{ord}_p(z)-l_pm_p)/2}\prod_{p \notin S(d)}p^{\mathrm{ord}_p(z)/2} \label{eq-form-plike} 
\end{equation} 
for some $m_p \in \mathbb N$ such that \eqref{eq-exp-cond}, $\xi _p^* \in \{\xi _p,\xi _p{}'\}$ $(p \in S(d)),$ and $m \in \mathbb Z$ such that 
\begin{equation} 
m \equiv \begin{cases} 
0,\ \pm 1 \pmod 3 & \text{if }\eta \in \mathbb Z[\sqrt d]\text{ or \ \:}z \equiv 0\ (\mathrm{mod}\ 2), \\ 
0\qquad\:\!\pmod 3 & \text{if }\eta \notin \mathbb Z[\sqrt d]\text{ and }z \equiv 1\ (\mathrm{mod}\ 2) 
\end{cases} \label{eq-trinity} 
\end{equation} 
and 
\begin{equation} 
1 \equiv \begin{cases} 
m\qquad\quad\:\:\:\pmod 2 & \text{if }x^2-dy^2 = -1\ \text{has an integral solution}, \\ 
\displaystyle\sum_ {p \in S(d)_-}m_p \pmod 2 & \text{otherwise}. 
\end{cases} \label{eq-cond-idx} 
\end{equation} 
\end{thm}
To prove Theorem \ref{thm-form-plike}, we use the following theorem.
\begin{thm}[{\cite[Theorem 2(a)]{Hir24}}]\label{thm-form-pell}
Let $d > 1$ be a square-free integer, and let $z > 1$ be an integer. 
Then $|x^2-dy^2| = z$ has a strictly primitive integral solution if and only if 
\[\mathrm{ord}_p(z) = \left\{\begin{array}{lll} 
l_pm_p & \text{if }p \in S(d), \\ 
0 & \text{if }p \notin S(d) 
\end{array}\right.\] 
for some $m_p \in \mathbb N$ for each prime number $p,$ where $l_2 = 2$ and $m_2 \in \{ 0,1\}$ if $\eta \notin \mathbb Z[\sqrt d].$ 
In this case, its strictly primitive integral solution $(x,y)$ satisfies 
\[ x+y\sqrt d = \pm\eta ^m\prod_{p \in S(d)}\xi _p^*{}^{m_p}\] 
for some $\xi _p^* \in \{\xi _p,\xi _p'\}$ $(p \in S(d))$ and $m \in \mathbb Z$ such that \eqref{eq-trinity}. 
\end{thm}
\begin{pf}[Proof of Theorem \ref{thm-form-plike}]
Let $(x,y) = (x_{\sigma (n)},x_{n+1})$ be an integral solution of $x^2-dy^2 = -z$ such that $\mathrm{gcd}\,(x_{\sigma (n)},d) = 1,$ and let $g = \mathrm{gcd}(x_{\sigma (n)},x_{n+1}).$ 
Then $(x,y) = (x_{\sigma (n)}/g,x_{n+1}/g)$ is a strictly primitive integral solution of $x^2-dy^2 = -z/g^2$ and satisfies 
\[\frac{x_{\sigma (n)}}{g}+\frac{x_{n+1}}{g}\sqrt d = \pm\eta ^m\prod_{p \in S(d)}\xi _p^*{}^{m_p}\] 
for some $m_p \in \mathbb N,$ $\xi _p^* \in \{\xi _p,\xi _p{}'\}$ $(p \in S(d)),$ and $m \in \mathbb Z$ such that \eqref{eq-trinity} by Theorem \ref{thm-form-pell}. 
Multiplying both sides by $g,$ we obtain 
\[ x_{\sigma (n)}+x_{n+1}\sqrt d = \pm\eta ^m\prod_{p \in S(d)}\xi _p^*{}^{m_p}p^{\mathrm{ord}_p(g)}\prod_{p \notin S(d)}p^{\mathrm{ord}_p(g)}.\] 
Taking norms of both sides, we obtain 
\[\mathrm{ord}_p(z) = \begin{cases} 
l_pm_p+2\,\mathrm{ord}_p(g) & \text{if }p \in S(d), \\ 
2\,\mathrm{ord}_p(g) & \text{if }p \notin S(d), 
\end{cases}\] 
which implies 
\[\mathrm{ord}_p(g) = \begin{cases} 
(\mathrm{ord}_p(z)-l_pm_p)/2 & \text{if }p \in S(d), \\ 
\mathrm{ord}_p(z)/2 & \text{if }p \notin S(d), 
\end{cases}\] 
\eqref{eq-exp-cond}, and \eqref{eq-form-plike}. 
Condition \eqref{eq-cond-idx} follows from 
\[\begin{cases} 
N(\xi _p^*) < 0 & \text{if }p \in S(d)_-, \\ 
N(\xi _p^*) > 0 & \text{if }p \in S(d)\setminus S(d)_- 
\end{cases}\] 
and 
\[\quad \begin{cases} 
N(\eta ) < 0 & \text{if }x^2-dy^2 = -1\ \text{has an integral solution}, \\ 
N(\eta ) > 0 & \text{otherwise}. 
\end{cases}\] 
\end{pf}
\begin{eg}\label{eg-sol-pell}
\begin{enumerate}
\item[(1)]
Every integral solution of $x_1{}^2+x_2{}^2+x_3{}^2 = 2x_4{}^2$ such that $(x_1,x_2) = (1,4)$ is given by 
\[ (x_1,x_2,x_3,x_4) = (1,4,\pm (5f_{2i+1}^{(2)}\pm 8g_{2i+1}^{(2)}),\pm (2f_{2i+1}^{(2)}\pm 5g_{2i+1}^{(2)})),\] 
where $i \in \mathbb N$ and the signs are chosen independently, since $x_3{}^2-2x_4{}^2 = -17,$ or equivalently, 
\[ x_3+x_4\sqrt 2 = \pm\eta ^{2i+1}\xi _{17}^*,\] 
where $\eta = 1+\sqrt 2,$ $\xi _{17} = 5+2\sqrt 2,$ and $\xi _{17}^* \in \{\xi _{17},\xi _{17}{}'\}.$
\item[(2)]
Every integral solution of $x_1{}^2+x_2{}^2+x_3{}^2 = 2x_4{}^2$ such that $(x_1,x_2) = (4,5)$ is given by 
\[ (x_1,x_2,x_3,x_4) = (4,5,\pm (7f_{2i+1}^{(2)}\pm 8g_{2i+1}^{(2)}),\pm (2f_{2i+1}^{(2)}\pm 7g_{2i+1}^{(2)})),\] 
where $i \in \mathbb N$ and the signs are chosen independently, since $x_3{}^2-2x_4{}^2 = -41,$ or equivalently, 
\[ x_3+x_4\sqrt 2 = \pm\eta ^{2i+1}\xi _{41}^*,\] 
where $\eta = 1+\sqrt 2,$ $\xi _{41} = 7+2\sqrt 2,$ and $\xi _{41}^* \in \{\xi _{41},\xi _{41}{}'\}.$
\end{enumerate}
\end{eg}
Combining Theorem \ref{thm-ang} with Theorem \ref{thm-form-plike}, we obtain the following solutions to Problem \ref{q-int-bisec}.
\begin{eg}
If $(x_1,x_2,x_3,x_4) = (a_1,a_2,a_3,a_4),$ $(b_1,b_2,b_3,b_4)$ are integral solutions of $x_1{}^2+x_2{}^2+x_3{}^2 = 2x_4{}^2$ for some $a_4,$ $b_4 \in \mathbb Z\setminus\{ 0\},$ then $\bm{a} = (a_1,a_2,a_3),$ $\bm{b} = (b_1,b_2,b_3),$ and 
\[\bm{c} = (b_4\bm{a}\pm a_4\bm{b})/g\] 
satisfy \eqref{eq-star}, where $g = \mathrm{gcd}(a_4,b_4).$ 
We summarize a few examples in the following table. 
\begin{center}\begin{tabular}{|c||c|c|c|c|c|c|c|c|c|} \hline 
$\bm{a}$ & $(1,7,0)$ & $(1,1,0)$ & $(23,89,0)$ & $(1,7,0)$ & $(1,1,0)$ & $(1,7,0)$ & $(1,1,4)$ & $(1,1,4)$ & $(1,4,9)$ \\ 
$a_4$ & $5$ & $1$ & $65$ & $5$ & $1$ & $5$ & $3$ & $3$ & $7$ \\ 
$\bm{b}$ & $(7,0,17)$ & $(1,0,7)$ & $(47,0,79)$ & $(1,1,4)$ & $(1,1,4)$ & $(3,4,5)$ & $(3,4,5)$ & $(7,7,8)$ & $(3,5,8)$ \\ 
$b_4$ & $13$ & $5$ & $65$ & $3$ & $3$ & $5$ & $5$ & $9$ & $7$ \\ 
$g$ & $1$ & $1$ & $65$ & $1$ & $1$ & $5$ & $1$ & $3$ & $7$ \\ \hline 
\end{tabular}\end{center}
\end{eg}
\begin{eg}
If 
\[ (a_{\sigma (1)},a_{\sigma (2)},a_{\sigma (3)}) = (1,4,5f_{2i+1}^{(2)}+8g_{2i+1}^{(2)}) \quad\text{and}\quad (b_{\tau (1)},b_{\tau (2)},b_{\tau (3)}) = (4,5,7f_{2j+1}^{(2)}+8g_{2j+1}^{(2)}),\] 
where $\sigma,$ $\tau \in \mathfrak S_3$ and $i,$ $j \in \mathbb N,$ then $\bm{a} = (a_1,a_2,a_3),$ $\bm{b} = (b_1,b_2,b_3),$ and 
\[\bm{c} = (2f_{2j+1}^{(2)}+7g_{2j+1}^{(2)})\bm{a}\pm (2f_{2i+1}^{(2)}+5g_{2i+1}^{(2)})\bm{b}\] 
satisfy \eqref{eq-star}. 
These vectors are produced by the integral solutions of $x_1{}^2+x_2{}^2+x_3{}^2 = 2x_4{}^2$ in Example \ref{eg-sol-pell}.
\end{eg}
\begin{eg}
If 
\[ (a_{\sigma (1)},a_{\sigma (2)},a_{\sigma (3)}) = (0,1,f_{2i+1}^{(2)}) \quad\text{and}\quad (b_{\tau (1)},b_{\tau (2)},b_{\tau (3)}) = (t^2-u^2,2tu,(t^2+u^2)f_{2j+1}^{(2)}),\] 
where $\sigma,$ $\tau \in \mathfrak S_3,$ $i,$ $j \in \mathbb N,$ and $(t,u) \in \mathbb Z^2\setminus\{ (0,0)\},$ then $\bm{a} = (a_1,a_2,a_3),$ $\bm{b} = (b_1,b_2,b_3),$ and 
\[\bm{c} = (t^2+u^2)g_{2j+1}^{(2)}\bm{a}\pm g_{2i+1}^{(2)}\bm{b}\] 
satisfy \eqref{eq-star}.
\end{eg}
\section{Incenters of Simplices over Fields}\label{sec-sim}
The following propositions are fundamental for simplices over $k.$ 
Let $O$ be the origin of $\mathbb R^n.$
\begin{prop}\label{prop-vol-sim}
In $\mathbb R^n,$ the volume $V$ of $n$-simplex over $k$ is $k$-rational.
\end{prop}
\begin{pf}
The volume $V$ of $n$-simplex $A_0A_1\cdots A_n$ over $k$ is given by 
\[ V = \frac{1}{n!}|\det (\overrightarrow{A_0A_1}\ \cdots\ \overrightarrow{A_0A_j}\ \cdots\ \overrightarrow{A_0A_n})|\] 
(see \cite{Ste66}), which implies $V \in k.$
\end{pf}
\begin{prop}\label{prop-rat-pt}
In $\mathbb R^n,$ with respect to $n$-simplex $A_0A_1\cdots A_n$ over $k,$ a point $P$ is $k$-rational if and only if $P$ has barycentric coordinates $(\lambda _0:\lambda _1:\cdots :\lambda _n)$ such that 
\begin{equation} 
\frac{\lambda _0}{\sum _{i = 0}^n\lambda _i},\ \frac{\lambda _1}{\sum _{i = 0}^n\lambda _i},\ \dots,\ \frac{\lambda _n}{\sum _{i = 0}^n\lambda _i} \in k. \label{eq-rat} 
\end{equation} 
In particular, if a point $P$ has barycentric coordinates $(\lambda _0:\lambda _1:\cdots :\lambda _n)$ such that $\lambda _0,$ $\lambda _1,$ $\dots,$ $\lambda _n \in k,$ then $P$ is $k$-rational.
\end{prop}
\begin{pf}
If $P \in k^n,$ then $P$ has barycentric coordinates $(\lambda _0:\lambda _1:\cdots :\lambda _n)$ such that \eqref{eq-rat} by letting $\lambda _i$ be the volume of the $n$-simplex with vertices $P$ and $A_j$ $(j \in \{ 0,1,\dots,n\}\setminus\{ i\})$ for each $i \in \{ 0,1,\dots,n\},$ where $\lambda _0,$ $\lambda _1,$ $\dots,$ $\lambda _n \in k$ by Proposition \ref{prop-vol-sim}.\par
Conversely, if $P$ has barycentric coordinates $(\lambda _0:\lambda _1:\cdots :\lambda _n)$ such that \eqref{eq-rat}, then $P \in k^n,$ since we have 
\[\overrightarrow{OP} = \frac{\sum_{i = 0}^n\lambda _i\overrightarrow{OA_i}}{\sum_{i = 0}^n\lambda _i}\] 
by definition of barycentric coordinates.\par
The second statement follows immediately from the first. 
\end{pf}
We are now ready to prove Theorem \ref{thm-sim}.
\begin{pf}[Proof of Theorem \ref{thm-sim}]
The incenter $I$ has barycentric coordinates 
\[ I = (a_0:a_1:\cdots :a_n).\]\par 
First, (S1) implies (S2) by Theorem \ref{thm-ang}, since $I$ is the intersection point of the bisectors of $\angle A_iA_hA_j$ $(h \neq i \neq j \neq h).$\par
Conversely, if (S2) holds, then there exist positive numbers $a_{0,0},$ $a_{1,0},$ $\dots,$ $a_{n,0} \in k$ such that $a_0 = a_{0,0}\sqrt d,$ $a_1 = a_{1,0}\sqrt d,$ $\dots,$ $a_n = a_{n,0}\sqrt d,$ which imply 
\[\frac{a_i}{a_0+a_1+\cdots +a_n} = \frac{a_{i,0}}{a_{0,0}+a_{1,0}+\cdots +a_{n,0}} \in k\] 
for each $i \in \{ 0,1,\dots,n\}$ and therefore $I \in k^n.$\par
The latter statement follows immediately from the former.
\end{pf}
\begin{eg}
The incenter of a given Heronian tetrahedron over $k$ is $k$-rational.
\end{eg}
\begin{eg}
Let $p,$ $q,$ $r \in k^\times.$ 
The equifacial tetrahedron with vertices $(p,q,r),$ $(p,-q,-r),$ $(-p,q,-r),$ and $(-p,-q,r)$ has a $k$-rational incenter at $(0,0,0).$
\end{eg}
Let $[XYZ]$ denote the area of triangle $XYZ.$
\begin{eg}
Let $p,$ $q \in k$ with $p \neq 3q$ and $p > 0.$ 
Consider the right regular triangular pyramid $ABCD$ with $A = (p,0,0),$ $B = (0,p,0),$ $C = (0,0,p),$ and $D = (q,q,q).$ 
Since $\overrightarrow{DB} = (q,q-p,q),$ $\overrightarrow{DC} = (q,q,q-p),$ and $\overrightarrow{DB}\times\overrightarrow{DC} = p(p-2q,q,q),$ we have $[BCD] = (1/2)|\overrightarrow{DB}\times\overrightarrow{DC}| = (p/2)\sqrt{(p-2q)^2+2q^2}.$ 
Similarly, we have $[CDA] = [DAB] = (p/2)\sqrt{(p-2q)^2+2q^2}.$ 
Furthermore, since triangle $ABC$ is the regular triangle with side length $p\sqrt 2,$ we have $[ABC] = (\sqrt 3/4)\cdot (p\sqrt 2)^2 = p^2\sqrt 3/2.$ 
Therefore, tetrahedron $ABCD$ has a $k$-rational incenter $I$ if and only if 
\[ (p-2q)^2+2q^2 \equiv 3 \pmod{k^{\times 2}}.\] 
If there exists a positive number $r \in k$ such that 
\[ (p-2q)^2+2q^2 = 3r^2,\] 
then we have 
\[ [BCD]:[CDA]:[DAB]:[ABC] = \frac{pr}{2}\sqrt 3:\frac{pr}{2}\sqrt 3:\frac{pr}{2}\sqrt 3:\frac{p^2}{2}\sqrt 3 = r:r:r:p\] 
and 
\[\overrightarrow{OI} = \frac{r(p,0,0)+r(0,p,0)+r(0,0,p)+p(q,q,q)}{r+r+r+p} = \frac{p(q+r)}{p+3r}(1,1,1),\] 
which implies $I$ is $k$-rational. 
For example, tetrahedron $ABCD$ with $A = (7,0,0),$ $B = (0,7,0),$ $C = (0,0,7),$ and $D = (1,1,1)$ has a rational incenter $I = (7/4,7/4,7/4),$ since $(7-2\cdot 1)^2+2\cdot 1^2 = 3\cdot 3^2.$
\end{eg}
\section{Centers of Triangles over Fields}\label{sec-trig}
In this section, we consider a few centers of triangle $ABC$ over $k$ with side lengths $BC = a,$ $CA = b,$ and $AB = c,$ perimeter $a+b+c = 2s,$ and area $\Delta.$ 
The following propositions are fundamental.
\begin{prop}\label{prop-isom}
In $\mathbb R^2,$ with respect to triangle $ABC$ over $k,$ the isotomic conjugate point $Q$ of a given $k$-rational point $P$ is also $k$-rational.
\end{prop}
\begin{pf}
Suppose that $P$ has normalized barycentric coordinates $(\lambda :\mu :\nu ).$ 
Then $Q$ has barycentric coordinates $(1/\lambda :1/\mu :1/\nu ).$ 
Thus, $P \in k^2$ implies that $\lambda,$ $\mu,$ $\nu \in k,$ and $1/\lambda,$ $1/\mu,$ $1/\nu \in k,$ and therefore $Q \in k^2.$
\end{pf}
\begin{prop}\label{prop-isog}
In $\mathbb R^2,$ with respect to triangle $ABC$ over $k,$ the isogonal conjugate point $Q$ of a given $k$-rational point $P$ is also $k$-rational.
\end{prop}
\begin{pf}
Suppose that $P$ has normalized barycentric coordinates $(\lambda :\mu :\nu ).$ 
Then $Q$ has barycentric coordinates $(a^2/\lambda :b^2/\mu :c^2/\nu ).$ 
Thus, $P \in k^2$ implies that $\lambda,$ $\mu,$ $\nu \in k,$ and $a^2/\lambda,$ $b^2/\mu,$ $c^2/\nu \in k,$ and therefore $Q \in k^2.$
\end{pf}
\begin{prop}
In $\mathbb R^2,$ the centroid $G,$ the Lemoine point $L,$ the circumcenter $E,$ and the orthocenter $H$ of triangle $ABC$ over $k$ are $k$-rational.
\end{prop}
\begin{pf}
As is well-known (see \cite{ETC}), the points $G,$ $L,$ $E,$ and $H$ have barycentric coordinates 
\begin{align*} 
G &= (1:1:1), \\ 
L &= (a^2:b^2:c^2), \\ 
E &= (a^2(b^2+c^2-a^2):b^2(c^2+a^2-b^2):c^2(a^2+b^2-c^2)), \\ 
H &= \left(\frac{1}{b^2+c^2-a^2}:\frac{1}{c^2+a^2-b^2}:\frac{1}{a^2+b^2-c^2}\right). 
\end{align*} 
Since $a^2,$ $b^2,$ $c^2 \in k,$ all of these coordinates lie in $k,$ which implies $G,$ $L,$ $E,$ $H \in k^2.$
\end{pf}
\begin{prop}
In $\mathbb R^2,$ the circumradius $R$ and the inradius $r$ of Heronian triangle $ABC$ over $k$ lie in $k.$
\end{prop}
\begin{pf}
If $a,$ $b,$ $c,$ $\Delta \in k,$ then $R,$ $r \in k$ since $\Delta = abc/4R = rs$ (see \cite[Section 1.5]{Cox89}).\qedhere
\end{pf}
We are now ready to prove Theorem \ref{thm-trig}.
\begin{pf}[Proof of Theorem \ref{thm-trig}]
As is well-known (see \cite{ETC}), the points $I,$ $I_A,$ $I_B,$ $I_C,$ and $Na$ have barycentric coordinates 
\begin{align*} 
I &= (a:b:c), \\ 
I_A &= (-a:b:c), \quad I_B = (a:-b:c), \quad I_C = (a:b:-c), \\ 
Na & = (s-a:s-b:s-c). 
\end{align*}\par
It follows from Theorem \ref{thm-sim} that $I \in k^2$ and \eqref{eq-trig} are equivalent.\par
We prove that $I \in k^2$ and $I_A \in k^2$ are equivalent. 
The condition $I \in k^2$ is equivalent to $a/2s,$ $b/2s,$ $c/2s \in k,$ and $-a/2(s-a),$ $b/2(s-a),$ $c/2(s-a) \in k,$ and therefore $I_A \in k^2.$ 
Here, we used the identities 
\[\left(\frac{a}{2s}\right) ^{-1}-\left(\frac{a}{2(s-a)}\right) ^{-1} = 2, \quad \left(\frac{b}{2s}\right) ^{-1}-\left(\frac{b}{2(s-a)}\right) ^{-1} = \frac{2a}{b}, \quad\text{and}\quad \left(\frac{c}{2s}\right) ^{-1}-\left(\frac{c}{2(s-a)}\right) ^{-1} = \frac{2a}{c}.\] 
Similarly, the condition $I \in k^2$ is equivalent to $I_B \in k^2,$ and $I_C \in k^2.$\par
We prove that $Na \in k^2$ implies $I \in k^2.$ 
If $Na \in k^2,$ then $(s-a)/s,$ $(s-b)/s,$ $(s-c)/s \in k$ since $(s-a)+(s-b)+(s-c) = s.$ 
This implies that $a/2s,$ $b/2s,$ $c/2s \in k,$ and $I \in k^2.$ 
The converse is also true.\par
Furthermore, $Ge \in k^2$ and $Na \in k^2$ are equivalent, since $Ge$ and $Na$ are isotomic to each other.\par
In addition, if \eqref{eq-trig} holds, then there exist positive numbers $a_0,$ $b_0,$ $c_0,$ $d \in k$ such that $a = a_0\sqrt d,$ $b = b_0\sqrt d,$ and $c = c_0\sqrt d,$ which imply $2s = (a_0+b_0+c_0)\sqrt d,$ $2(s-a) = (b_0+c_0-a_0)\sqrt d,$ $2(s-b) = (c_0+a_0-b_0)\sqrt d,$ $2(s-c) = (a_0+b_0-c_0)\sqrt d,$ and therefore $r^2,$ $r_A{}^2,$ $r_B{}^2,$ and $r_C{}^2$ are numbers in $k^\times$ equivalent to $d$ modulo $k^{\times 2},$ since $\Delta \in k$ and $\Delta = rs = r_A(s-a) = r_B(s-b) = r_C(s-c)$ (see \cite[Section 1.5]{Cox89}).  
This implies that triangle $ABC$ is similar to the Heronian triangle over $k$ with side lengths $a_0,$ $b_0,$ and $c_0,$ and area $\Delta /d.$
\end{pf}
From Theorem \ref{thm-trig}, we obtain the following construction method for every triangle over $k$ with $k$-rational incenter.
\begin{cor}
In $\mathbb R^2,$ every triangle $ABC$ over $k$ with $k$-rational incenter can be constructed by the following steps.\par
{\tt Step 1.} 
Choose a positive number $d \in k.$\par
{\tt Step 2.} 
Choose solutions $(x_1,x_2,x_3) = (u_1,u_2,u_3),$ $(v_1,v_2,v_3),$ $(w_1,w_2,w_3)$ of $x_1{}^2+x_2{}^2 = dx_3{}^2$ in $k$ such that $\bm{u} = (u_1,u_2),$ $\bm{v} = (v_1,v_2),$ and $\bm{w} = (w_1,w_2)$ are pairwise linearly independent over $k.$\par
{\tt Step 3.} 
Find the intersection point $C_0$ of the line passing through $B_0 = (w_1,w_2)$ and parallel to $\bm{u}$ and the line connecting $O$ and $(v_1,v_2).$\par
{\tt Step 4.} 
Translate triangle $OB_0C_0$ by a certain vector in $k^2.$
\end{cor}
\begin{pf}
Consider triangle $ABC$ over $k$ with vertices $A = (a_1,a_2),$ $B = (b_1,b_2),$ and $C = (c_1,c_2),$ and $k$-rational incenter. 
Since 
\[ |\overrightarrow{BC}|^2 \equiv |\overrightarrow{CA}|^2 \equiv |\overrightarrow{AB}|^2 \pmod{k^{\times 2}}\] 
by Theorem \ref{thm-trig}, there exist numbers $u_3,$ $v_3,$ $w_3,$ $d \in k$ such that $d > 0$ and $(x_1,x_2,x_3) = (c_1-b_1,c_2-b_2,u_3),$ $(a_1-c_1,a_2-c_2,v_3),$ $(b_1-a_1,b_2-a_2,w_3)$ are solutions of $x_1{}^2+x_2{}^2 = dx_3{}^2$ in $k.$ 
Therefore, let $\bm{u} = \overrightarrow{BC},$ $\bm{v} = \overrightarrow{CA},$ and $\bm{w} = \overrightarrow{AB},$ which are pairwise linearly independent over $k.$ 
The line passing through $B_0 = (b_1-a_1,b_2-a_2)$ and parallel to $\bm{u}$ and the line connecting $O$ and $(a_1-c_1,a_2-c_2)$ intersect at $C_0 = (c_1-a_1,c_2-a_2).$ 
Triangle $ABC$ can be obtained by translating triangle $OB_0C_0$ by vector $\overrightarrow{OA} = (a_1,a_2).$
\end{pf}
\begin{eg}
In the following cases, the centroid $G,$ the Lemoine point $L,$ the circumcenter $E,$ the orthocenter $H,$ the incenter $I,$ the excenters $I_A,$ $I_B,$ and $I_C,$ the Gergonne point $Ge,$ and the Nagel point $Na$ of triangle $ABC$ are $k$-rational. 
For readers' convenience, we add the values of the vectors $\bm{u} = \overrightarrow{BC},$ $\bm{v} = \overrightarrow{CA},$ and $\bm{w} = \overrightarrow{AB},$ the side lengths $a,$ $b,$ and $c,$ the inradius $r,$ the exradii $r_A,$ $r_B,$ and $r_C,$ the positive algebraic integer $d \in O_k$ such that $a^2 \equiv b^2 \equiv c^2 \equiv d \pmod{k^{\times 2}},$ and the area $\Delta$ (see Figure 2 for the example in the middle column).
\begin{center}\begin{tabular}{|c||c|c|c|} \hline 
$k$ & $\mathbb Q$ & $\mathbb Q$ & $\mathbb Q(\sqrt 2)$ \\ \hline 
$\bm{u}$ & $5(-3,4)$ & $14(-1,1)$ & $7(-1,\sqrt 2)$ \\ 
$\bm{v}$ & $-3(5,12)$ & $-3(1,7)$ & $-10(5,\sqrt 2)$ \\ 
$\bm{w}$ & $2(15,8)$ & $(17,7)$ & $3(19,\sqrt 2)$ \\ \hline 
$A$ & $(0,0)$ & $(0,0)$ & $(0,0)$ \\ 
$B$ & $(30,16)$ & $(17,7)$ & $(57,3\sqrt 2)$ \\ 
$C$ & $(15,36)$ & $(3,21)$ & $(50,10\sqrt 2)$ \\ \hline 
$G$ & $(15,52/3)$ & $(20/3,28/3)$ & $(107/3,13\sqrt 2/3)$ \\ 
$L$ & $(31485/1651,32976/1651)$ & $(2166/295,2562/295)$ & $(52875/1019,6795\sqrt 2/1019)$ \\ 
$E$ & $(72/7,943/56)$ & $(47/8,79/8)$ & $(117/4,-45\sqrt 2/8)$ \\ 
$H$ & $(171/7,513/28)$ & $(33/4,33/4)$ & $(97/2,97\sqrt 2/4)$ \\ \hline 
$I$ & $(120/7,132/7)$ & $(7,9)$ & $(48,6\sqrt 2)$ \\ 
$I_A$ & $(35,77/2)$ & $(21,27)$ & $(60,15\sqrt 2/2)$ \\ 
$I_B$ & $(-33,30)$ & $(-18,14)$ & $(-6,24\sqrt 2)$ \\ 
$I_C$ & $(22,-20)$ & $(27/2,-21/2)$ & $(15,-60\sqrt 2)$ \\ 
$Ge$ & $(96/5,96/5)$ & $(539/73,637/73)$ & $(5096/103,784\sqrt 2/103)$ \\ 
$Na$ & $(75/7,100/7)$ & $(6,10)$ & $(11,\sqrt 2)$ \\ \hline 
$a$ & $25$ & $14\sqrt 2$ & $7\sqrt 3$ \\ 
$b$ & $39$ & $15\sqrt 2$ & $33\sqrt 3$ \\ 
$c$ & $34$ & $13\sqrt 2$ & $30\sqrt 3$ \\ 
$r$ & $60/7$ & $4\sqrt 2$ & $2\sqrt 6$ \\ 
$r_A$ & $35/2$ & $12\sqrt 2$ & $5\sqrt 6/2$ \\ 
$r_B$ & $42$ & $14\sqrt 2$ & $14\sqrt 6$ \\ 
$r_C$ & $28$ & $21\sqrt 2/2$ & $35\sqrt 6$ \\ \hline 
$d$ & $1$ & $2$ & $3$ \\ \hline 
$\Delta$ & $420$ & $168$ & $210\sqrt 2$ \\ \hline 
\end{tabular}\end{center}
\begin{figure}[h]
\centering
\includegraphics{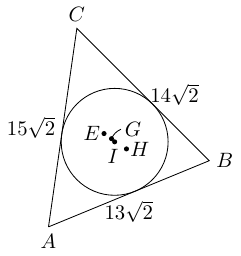}
\caption{The centroid $G,$ the circumcenter $E,$ the orthocenter $H,$ and the incenter $I$ of triangle $ABC$ with $A = (0,0),$ $B = (17,7),$ and $C = (3,21).$}
\end{figure}
\end{eg}


\begin{thebibliography}{99}
\bibitem{Car15}
R.~D.~Carmichael, {\it Diophantine Analysis}, John Wiley \& Sons, New York, 1915.
\bibitem{Cox89}
H.~S.~M.~Coxeter, {\it Introduction to Geometry}, 2nd ed., John Wiley \& Sons, New York, 1989.
\bibitem{Hir24}
T.~Hirotsu, General Pell's equations and angle bisectors between planar lines with rational slopes, {\it Integers}, {\bf 24} (2024), \#A111.
\bibitem{ETC}
C.~Kimberling, Encyclopedia of Triangle Centers, \\
{\tt https://faculty.evansville.edu/ck6/encyclopedia/ETC.html} (retrieved 15 Dec.~2025).
\bibitem{MP13}
S.~H.~Marshall and A.~R.~Perlis, Heronian tetrahedra are lattice tetrahedra, Amer.~Math.~Monthly, {\bf 120} (2), 2013, 140--149.
\bibitem{Som58}
D.~M.~Y.~Sommerville, {\it An Introduction to the Geometry of $n$ Dimensions}, Dover, New York, 1958.
\bibitem{Ste66}
P.~Stein, A note on the volume of a simplex, Amer.~Math.~Monthly, {\bf 73} (3), 1966, 299--301.
\end{thebibliography}
\end{document}